\documentclass{ifacconf}

\usepackage{graphicx}      
\usepackage{natbib}        

\usepackage{graphics} 
\usepackage{epsfig} 
\usepackage{algorithm} 
\usepackage{algpseudocode} 
 
\usepackage{amsmath}
\usepackage{amssymb}
\usepackage{gensymb}
\usepackage{optidef}

\DeclareMathOperator*{\argmin}{argmin}

\begin{document}
\begin{frontmatter}

\title{Safe Optimization of an Industrial Refrigeration Process Using an Adaptive and Explorative Framework} 

\author[First]{Buse Sibel Korkmaz}
\author[Second]{Marta Zagórowska}
\author[First]{Mehmet Mercangöz}

\address[First]{Imperial College London, 
   London SW7 2AZ, UK (e-mail: \textit{\{buse.korkmaz18, m.mercangoz\}}@ imperial.ac.uk)}
\address[Second]{ETH Zürich, 8092 Zürich, Switzerland  (e-mail: mzagorowska@control.ee.ethz.ch)}

\begin{abstract} Many industrial applications rely on real-time optimization to improve key performance indicators. In the case of unknown process characteristics, real-time optimization becomes challenging, particularly for the satisfaction of safety constraints. In this paper, we demonstrate the application of an adaptive and explorative real-time optimization framework to an industrial refrigeration process, where we learn the process characteristics through changes in process control targets and through exploration to satisfy safety constraints. We quantify the uncertainty in unknown compressor characteristics of the refrigeration plant by using Gaussian processes and incorporate this uncertainty into the objective function of the real-time optimization problem as a weighted cost term. We adaptively control the weight of this term to drive exploration. The results of our simulation experiments indicate the proposed approach can help to increase the energy efficiency of the considered refrigeration process, closely approximating the performance of a solution that has complete information about the compressor performance characteristics.

\end{abstract}


\end{frontmatter}

\section{INTRODUCTION}

Many industrial applications use real-time optimization based on characteristics of the underlying processes to improve key performance indicators, such as energy consumption. However, the characteristics may be unknown, or may change with time, for instance due to wear and tear. Thus, model-based real-time optimization is challenging. A common approach is to first estimate the unknown characteristics by collecting observations from multiple operating points. The operating points are explored, for example, by doing plant step tests. These observations generally bring noise to the optimization problem, and a real-time optimization algorithm needs to sequentially estimate and optimize the utility function by utilizing these noisy observations to reveal the decisions leading to the lowest cost. 

Considering the cost of exploration such as not delivering the expected demand from the system during exploration, there is a trade-off between exploring more and having more accurate estimations or exploiting collected observations and optimizing the obtained function. This trade-off is well-studied in the literature and multi-armed bandit approaches are suggested as a solution for optimally balancing exploration and exploitation \citep{bubeck2012regret}. The first appearance of using upper confidence bounds (UCB) is dated back to \citet{1985Lai}. After that, many UCB algorithms are developed by researchers for stochastic bandits with many arms \citep{banditbook}. Most efficient bandit algorithms are characterized by having a utility function under certain regularity conditions \citep{Dani2008}. \citet{Srinavas2010} divided the stochastic optimization problem into two objectives as unknown function estimations from noisy observations and optimization of the estimated function over the decision set. To satisfy regularity conditions, they leverage the kernel methods and Gaussian process (GP) to model utility function through the regularity assumptions of kernels \citep{RasmussenGP}. 

If we want to approximate the unknown characteristics and optimize a process in a safety-critical operation, we need to utilize safe exploration algorithms which allow only the exploration of feasible decision points by enforcing safety constraints. In this work, we define a feasible decision point in the given decision set as satisfying all the constraints of the given problem. Bandit setting and Markov Decision Processes (MDPs) are used to formalize the safe exploration problem in the literature \citep{Schreiter2015, Turchetta2019}. \citet{Sui2015} have introduced the SafeOpt algorithm which optimizes an unknown function over a safe decision set notion under certain assumptions by utilizing the confidence bound construction method of \citet{Srinavas2010}. Following a similar approach we recently proposed the ARTEO algorithm \citep{ARTEO} to cast multi-armed bandits as a mathematical programming problem subject to safety constraints for the optimization of safety-critical systems. In ARTEO, we model the unknown characteristics of the system using Gaussian processes (GP) and utilize the covariance function of GPs to quantify uncertainty to ensure satisfaction of the safety constraints with high probability. We establish the safety of ARTEO by constructing  confidence bounds as \citet{Srinavas2010}. We incorporate the uncertainty into the utility function as a contribution to encourage the exploration at points with high uncertainty while continuing to satisfy the optimization goals of the system. The size of this contribution is adaptively controlled by a hyperparameter in accordance with the requirements of the environment.

In this paper, we demonstrate the implementation of ARTEO for the power management problem of an industrial refrigeration process.  The results show that we are able to successfully control three key performance indicators of the process by applying ARTEO. We further develop different demand scenarios to investigate the potential and limitations of ARTEO such as the case when the safety condition is violated, how hyperparameter selection affects the performance of the algorithm, and how ARTEO performs compared to a benchmark optimization problem with fully known system characteristics.

\section{PROBLEM STATEMENT AND BACKGROUND}

We want to minimize a utility function $f$ that models the cost of decisions and the uncertainty in the system subject to $K$ safety constraints $g_k$, $k=1,\ldots, K$. The cost of a decision and the constraints at the chosen point $x_t \in D$ at each solution iteration $t$ depends on the unknown system characteristics $p_t \in P$, and $D_i \subset \mathbb{R}$ defines the decision set of arm $i$ . We want to find a sequence of decisions, $x_1,x_2, \ldots, x_T$ which minimizes $f$ subject to $g$ while revealing the unknown characteristics. After making a decision $x_t$, we obtain a noisy measurement, $y_t=f(x_t,p(x_t))+\epsilon$ where $f:D \times P \rightarrow \mathbb{R}$ and $\epsilon$ is $R$-sub-Gaussian noise for a fixed constant $R \geq 0$ \citep{SubGaussian}. At every iteration, the decision $x_t$ must satisfy the  constraints $g_k(x_t,p(x_t))+ h_k\leq 0$, where $k=1,2,\ldots, K$ with $K$ denoting the number of constraints. The value of $h_k\geq 0$ is called a \emph{safety threshold}. Thus, we can formalise our optimisation problem at time $t$ as:

\begin{subequations}
\begin{align}
X_t=&\argmin_{X} \quad f(X,p(X))\\
\text{subject to} & \quad g_k(X,p(X))+h_k\leq 0, \forall k=1,\ldots,K
\end{align} 
\end{subequations}
where $X = [x_{i}, ..., x_{N}]$ for arms ${i=1, \dots, N}$. Since the characteristics $p(.)$ are unknown and need to be estimated, at every solution instance $t$, we need to first solve an estimation problem to find $p$, then solve the optimization problem (1). In ARTEO, we utilize GPs to learn $p(.)$ and its covariance matrix to establish confidence bounds that include the true value of safety function with $1 - \delta $ probability where $\delta \in (0,1)$ and a coefficient $\beta_t$ is chosen to satisfy explained conditions in \citet{ARTEO}. This guarantee can be formalized as next:

\begin{equation}
P\left[\left|g(x_t) - \mu_{t-1}(x_t)\right| \leq  \sqrt{\beta_t}\sigma_{t-1}(x_t)\right] \leq 1 - \delta, \ t \geq 1
\end{equation}

In \citet{ARTEO}, we established the theoretical guarantee of the confidence bounds by choosing $\beta_t$ based on Theorem 3 of \citet{Srinavas2010}, Theorem 2 of \citet{BoundProof}, and Theorem 1 of \citet{Sui2015} under the regularity assumptions of GPs. Further details regarding $\beta$ parameter and regularity assumptions can be found in \citet{ARTEO}. The implementation details of the algorithm are given in the next section. 

\subsection{Gaussian Processes}

Gaussian processes (GPs) are non-parametric models which can be used for regression. GPs are fully specified by a mean function $\mu(x)$ and a kernel  $k(x,x')$ which is a covariance function and decides the shape of prior and posterior in GP regression \citep{RasmussenGP}. The goal is to predict the value of $p$ including unknown components at the decision point $x^*$ by using GPs and then compute $f$. Assuming having a zero mean prior, the posterior $p(f(x^*)|x^*, D)$ follows $N(\mu(x^*), \sigma^2(x^*))$ that satisfy, 
\begin{equation}
\begin{aligned}
  \mu(x^*) = K_{N*}^T(K_{NN} + \sigma^2 I)^{-1}(y_1,...,y_N)^T\\
  \sigma^2(x^*) = k(x^*,x^*) - K_{N*}^T(K_{NN} + \sigma^2 I)^{-1}
\end{aligned}
\end{equation}

where $t,j \in \left\{1,.., N\right\}$ and denotes the index of observations, $[K_{N*}]_t = k(x_t,x^*)$, $K_{NN}$ is the $N\times N$ positive definite kernel matrix with  $[K_{NN}]_{t,j} = k(x_t, x_j)$ and $y_t$ is the noisy feedback of $t$th observation.

\section{ARTEO ALGORITHM}

ARTEO is developed for safety-critical systems with high exploration costs. It does not require a separate training phase and learns during normal operation. The algorithm is executed periodically and the optimal solutions for (1) can change due to a change in optimization targets, a change in system characteristics, or can be driven by the exploratory contribution in the cost function. At each solution instance ARTEO updates the posterior distributions of GPs with previous noisy observations and provides the optimal solution for the desired outcome based on how GPs model the unknown components. In a setup with multiple-arms, the decision set $D_i$ is defined for each arm $i$ as satisfying the introduced assumptions in \citet{ARTEO}. For each arm, a GP prior and initial \emph{safe seed} set is introduced to the algorithm. The safe seed set $S_0$ includes at least one safe decision point with the true value of the safety function at that point. Safe seed points are given to accurately assess the safety of chosen operating points in later iterations as in previously published safe exploration algorithms \citep{Sui2015, Sui2018,Turchetta2019}. At the first solution instance, each GP corresponding to arm $i$ is fit by a safe seed set and given GP prior to computing a posterior distribution. At later instances, the safe set is expanded by noisy observations as collected after the execution of decided operating points from the real-time optimizer. These observations are used to compute posteriors to use in the optimization of the utility function, which includes the cost of decision and uncertainty.  

The uncertainty in the utility function $f_t$ is quantified as:
\begin{equation}
\begin{array}{ll}
     &  U_t(X_t, p(X_t)) = \sum_{i}^{n}{\sigma(x_{it})}
\end{array}
\end{equation}
where $\sigma(x_{it})$ is the standard deviation of $GP_i$ at the point $x_{it}$ for the arm $i$ in the iteration $t$ and $n$ is the number of arms. It is incorporated into the utility function $f$ by multiplying by an adjustable parameter $z$ as next:
\begin{equation}
\begin{array}{ll}
    f_t(X_t, p(X_t)) = C_t(X_t, p(X_t)) + z U_t(X_t, p(X_t))
\end{array}
\end{equation}

where $C_t(X_t, p(X_t))$ represents the cost of decision at the evaluated points. The uncertainty weight $z$ is adapted using present uncertainty in the environment. The impact of adjusting $z$ on the performance of ARTEO will be shown in Section \ref{sec:HyperparameterSelection}. 


The optimization phase occurs after updating the posteriors of each arm. The real-time optimizer incorporates the posterior of the GP into decision-making by modelling the unknown parameters for each arm by using the mean and standard deviation of GPs. The utility function is the objective function in the real-time optimization formulation and the safety thresholds are constraints. In the utility function, the mean of the GP posterior of each arm is used to evaluate the cost of the decision and the standard deviation of GPs is used to measure uncertainty. In the safety function, the standard deviation of the GP posterior of each arm is used to construct confidence bounds and then these bounds are used to assess the feasibility of evaluated points. The optimizer solves the minimization problem under safety constraints within the defined decision set of each arm. Any optimization algorithm that could solve the given problem can be used in this phase. We have a constrained nonlinear case study, and we choose an interior-point algorithm to solve our problem.

Since we apply the ARTEO algorithm to a dynamical system simulation, we build a steady-state detector and evaluate costs and constraints in a steady-state. For the steady-state detection algorithm, we use a statistical test-based detector as common in many real-time optimization tools \citep{SteadyStateRTO}. The pseudocode of ARTEO is given in \citet{ARTEO}.

\begin{figure}[b]
  \centering
  \includegraphics[scale=0.35]{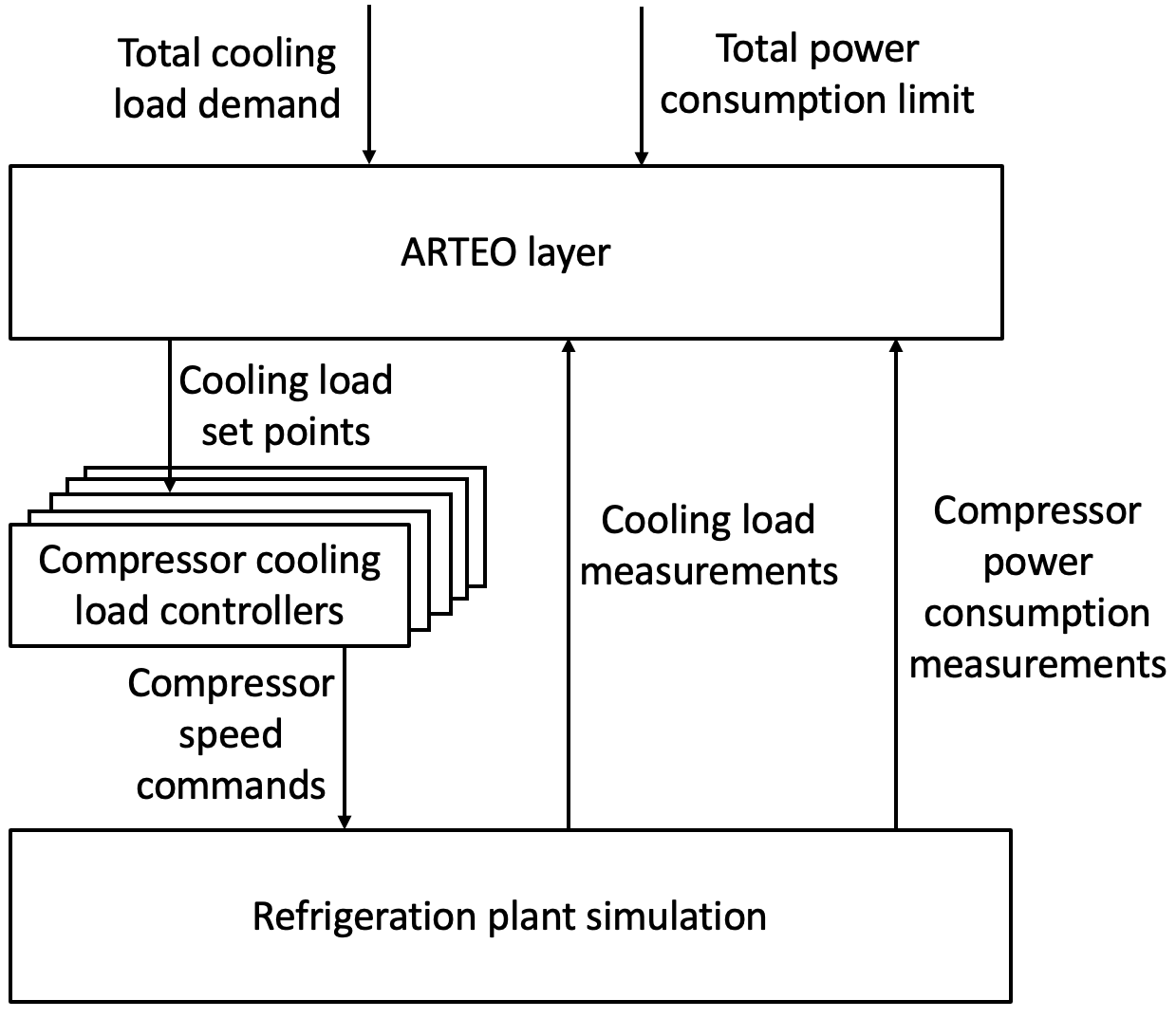}
  \caption{Simulation of the refrigeration system with ARTEO}
  \label{Figure1}
\end{figure}

\section{CASE STUDY}

We apply the ARTEO approach to a food freezing process plant. This case study is introduced first in \citet{fishfreezingpaper}, which we develop further into a dynamic simulation. We develop the simulation and algorithm with MATLAB/Simulink 2022a, and conduct experiments on an M1 Pro chip with 16 GB memory. The general structure of the industrial refrigeration control system with the ARTEO layer used in the simulations is given in Fig. \ref{Figure1}. 

\subsection{The industrial refrigeration system}

\begin{figure}[t]
  \centering
  \includegraphics[scale=0.135]{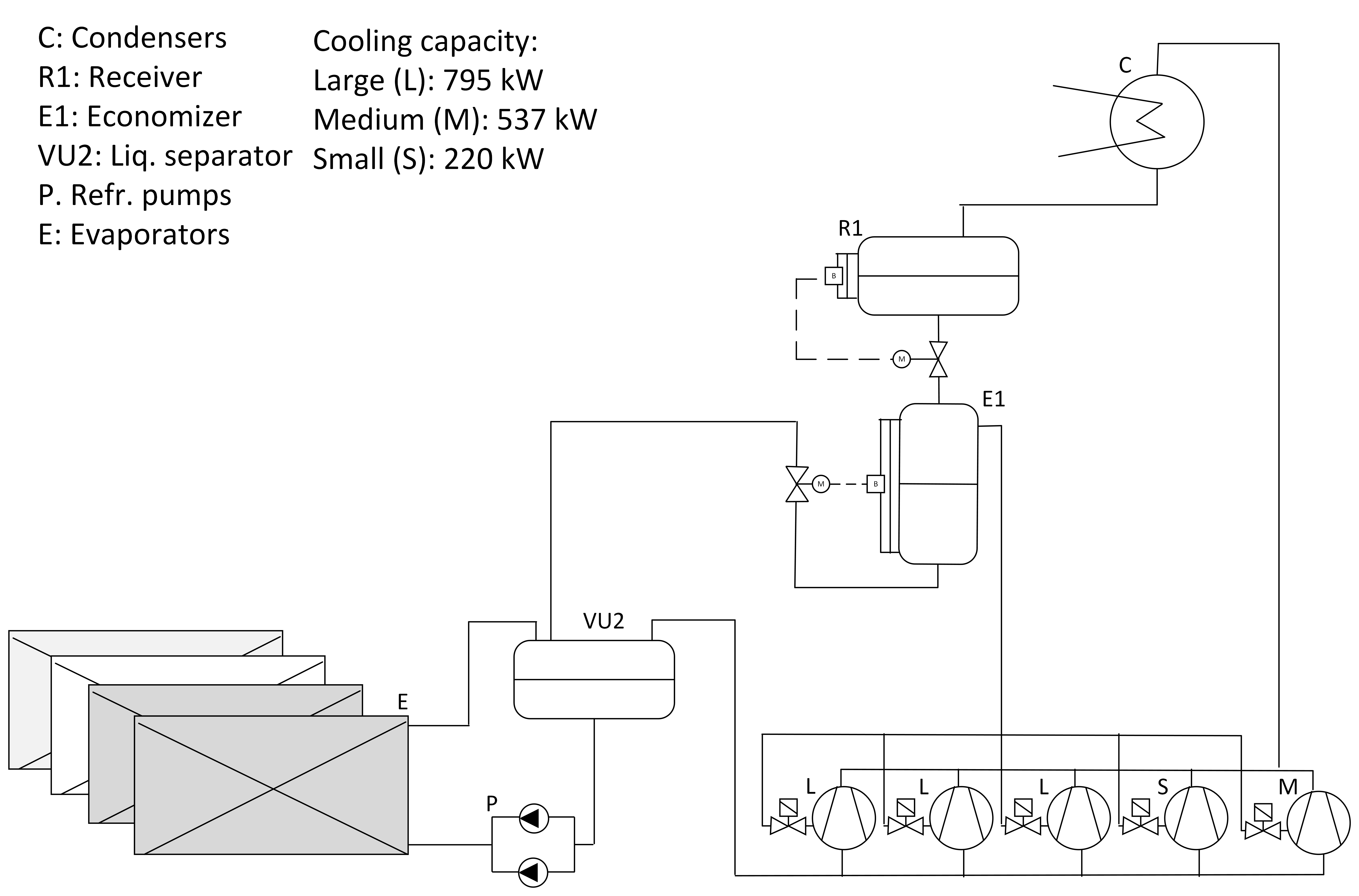}
  \caption{Schematic of the case study refrigeration system}
  \label{Figure2}
\end{figure}

\begin{figure}[b]
  \centering
  \includegraphics[scale=0.5]{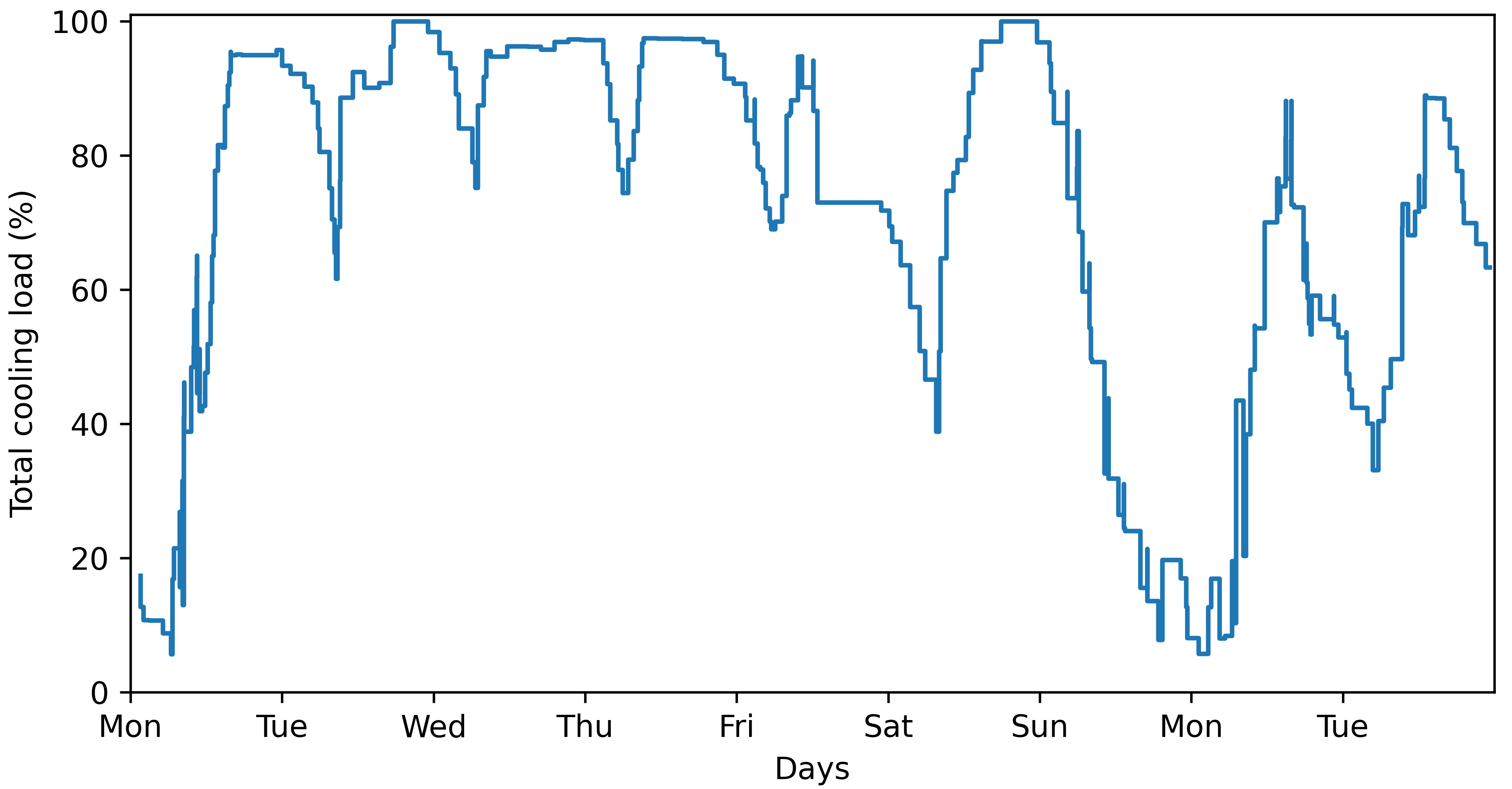}
  \caption{Cooling load variations in the case study process}
  \label{Figure3}
\end{figure}

\begin{figure*}[t]
  \centering
  \includegraphics[scale=0.34]{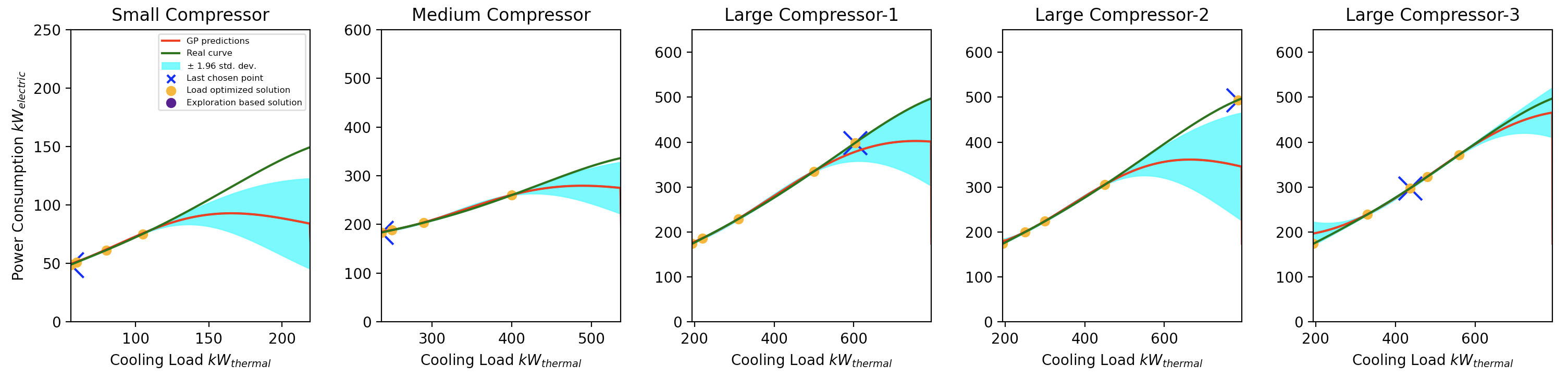}
  \caption{The first iteration of the algorithm. The safe seed set includes three points for each compressor within their respective operating intervals. The blue-shaded area represents the uncertainty which is high due to unknown regions}
  \label{Figure4}
\end{figure*}

\begin{figure*}[h]
  \centering
  \includegraphics[scale=0.34]{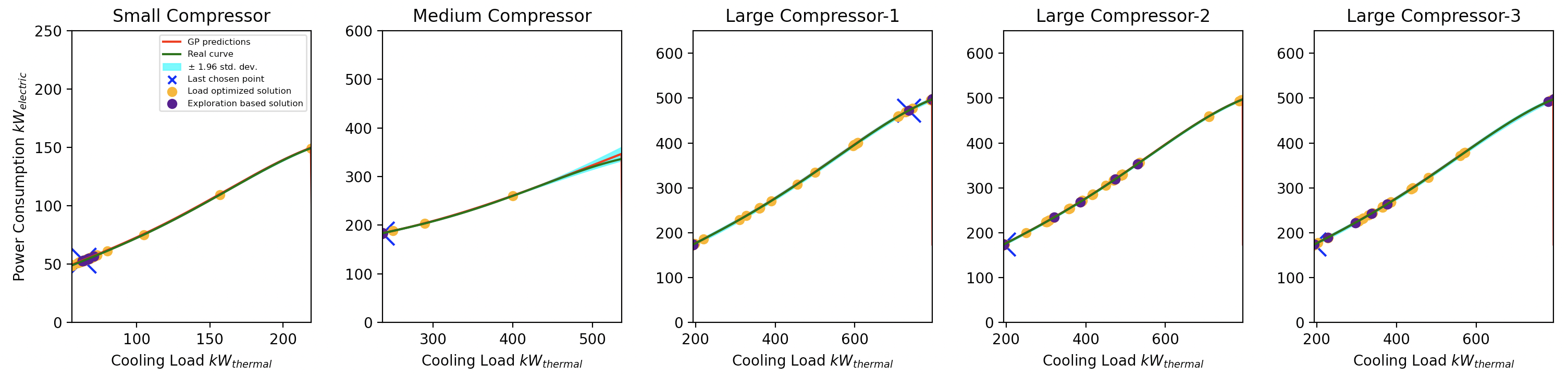}
  \caption{The final iteration in the simulation. GP predicted performance lines are converged to the actual performance curves of each compressor. Purple-coloured sample points are chosen by exploration when the environment becomes available}
  \label{Figure5}
\end{figure*}

The industrial refrigeration system is formed by two subsystems with different pressure levels. The two subsystems consist of 11 compressors: three reciprocating, and eight screw compressors used to freeze stored fish within the refrigeration system. We focus on the screw compressors because they have varying ranges of cooling capacities and different power consumption characteristics. The first subsystem has 5 screw compressors and the second subsystem has 3 screw compressors. The system with five screw compressors is chosen to apply our approach. Fig. \ref{Figure2} illustrates a simplified process layout of this subsystem.

For the dynamic modelling of the screw compressors, 200 seconds is considered to be the time needed to move from one steady-state to a new one based on the dynamic experiments from \citet{dynamiccompressor}. Therefore, ARTEO is triggered by 250 seconds time intervals to decouple the optimization actions from the compressor controllers. For representing the performance of screw compressors, the characteristics capturing the cooling load and power consumption data is obtained from \citet{perfmulcompressor} and embedded into simulations but not directly revealed to the optimization algorithm. The system consists of one small-sized, one medium-sized and three large-sized compressors. The minimum and maximum cooling capacities for each compressor type are given in Table \ref{Table1}. The coefficient of performance for cooling for all compressors is assumed to be 1.6 based on data from \cite{fishfreezingpaper}. 

The cooling demand is subject to change based on the stored amount of fish and environmental factors. Fig. \ref{Figure3} demonstrates the production load during fishing season when the plant is needed to be cooled at full capacity for most days but also with considerable variations. This period is chosen to highlight that large production rates require higher energy consumption, which can be a safety constraint in some cases due to a shortage of electricity in the power grid or similarly a critical cost concern when regulations or power purchase contracts demand enforcing a peak electricity consumption limit. 

More specifically, we demonstrate our approach over the period Saturday, Sunday and Monday in Fig. \ref{Figure3} because those days have significant fluctuations in the process cooling demand, which drive the plant to work at different operating points. Operating at unknown points may cause the violation of constraints due to unknown performance characteristics at those points.

\begin{table}[h]
\caption{Compressors Operating Intervals}
\label{Table1}
\begin{center}
\begin{tabular}{|p{0.20\linewidth}||p{0.32\linewidth}|p{0.32\linewidth}|}
\hline
Compressor & Min. Capacity (kW) & Max. Capacity (kW)\\
\hline
Small & 56 & 220\\
\hline
Medium & 237 & 537\\
\hline
Large & 194 & 795\\
\hline
\end{tabular}
\end{center}
\end{table}

\subsection{Safe learning and optimization with ARTEO}

For this case study, we choose the Squared Exponential kernel, which is a positive definite kernel function, as the GP prior for each compressor based on the suggestion of the GP hyperparameter optimization tool under the Statistics and Machine Learning Toolbox in Matlab \citep{RasmussenGP}. The minimum and maximum cooling capacities in Table 1 are introduced as decision-set boundaries to the optimization: 
\begin{equation}
    x_{i\min}\leq x_{it} \leq x_{i\max} \quad i = 1, \ldots, 5, t = 1, \ldots, 42000
\end{equation}
Then, we create the safe seed set including three safe decision points for each compressor to obtain the posterior of GPs before making any decisions. We define the utility function $f$ as follows:
\begin{equation}
f_t(x) = \left[\sum_{i=1}^{5}{\mu_{P_i(x_{it})}}\right]^2 + \left[M_t - \sum_{i=1}^{5}{x_{it}}\right]^2 - z\sum_{i=1}^{5}{\sigma_{P_i(x_{it})}}
\end{equation}
where $\mu_{P_i(x_{it})}$ represents the mean prediction of power consumption at production load $x$ and $\sigma_{P_i(x_{it})}$ represents the standard deviation of prediction to represent uncertainty in predictions at load $x_{it}$. $M_t$ denotes the desired production load at time $t$ based on the introduced cooling scenario. We develop the following two conditions and assign a positive value $z$ = 1000 when these conditions hold: (1) the cooling demand has stayed the same in the two previous sampling instances, (2) there is at least one feasible decision point that satisfies the cooling demand with a tolerance $\alpha$ = 10 kW. The safety threshold for total power consumption is set as 80\% of the maximum power consumption possible for all compressors, which corresponds to 1580 kW, to emulate an external command given to limit the power consumption of the process. The safety constraint $g$ is defined next as:
\begin{equation}
g(x) = \sum_{i=1}^{5}{\left[\mu_{P_i(x_{it})} + \sqrt{\beta_t}\sigma_{P_i(x_{it})}\right]}
\end{equation}
where $\beta_t$ is chosen as explained in \citet{ARTEO}.
\begin{figure}[t!]
  \centering
  \includegraphics[width=8.4cm]{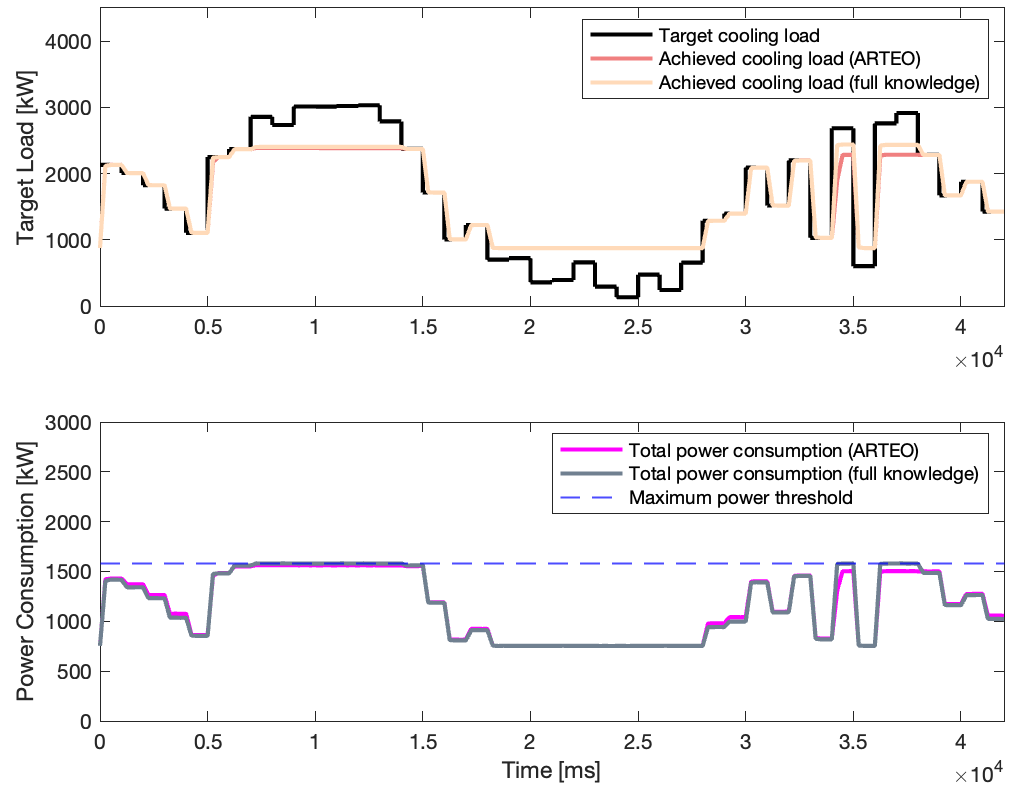}
  \caption{Target, achieved cooling demands and total power consumption over time}
  \label{Figure5-6}
\end{figure}
The experiment is simulated with the given cooling load trajectory in Fig. \ref{Figure2}. Fig. \ref{Figure3} demonstrates the first solution instance of the algorithm, where the uncertainty is high and it can be recognized from the wide confidence bounds. At each instance, GPs are updated with new observations, and the algorithm takes new decisions by the updated GPs. Fig. \ref{Figure4} shows the final results for performance curve estimations at the end of the simulated cooling load scenario. It can be seen that the ARTEO algorithm is able to successfully approximate the power consumption curves very close to their real shape. Fig. \ref{Figure5-6} illustrates the expected and achieved production loads. The expected production load is satisfied unless: (1) it is higher than what compressors can achieve without crossing the safety threshold, (2) it requires compressors to operate lower than their minimum operating points. Fig. \ref{Figure5-6} demonstrates that when the total power consumption is on the limit of the maximum power threshold, the expected load could not be satisfied. At those production loads, our approach finds safe operating points which minimize the stated objective function. 

\subsection{Limitations of confidence-bound-based safety}

\begin{figure}[tbhp]
    \begin{center}
    \includegraphics[width=8.4cm]{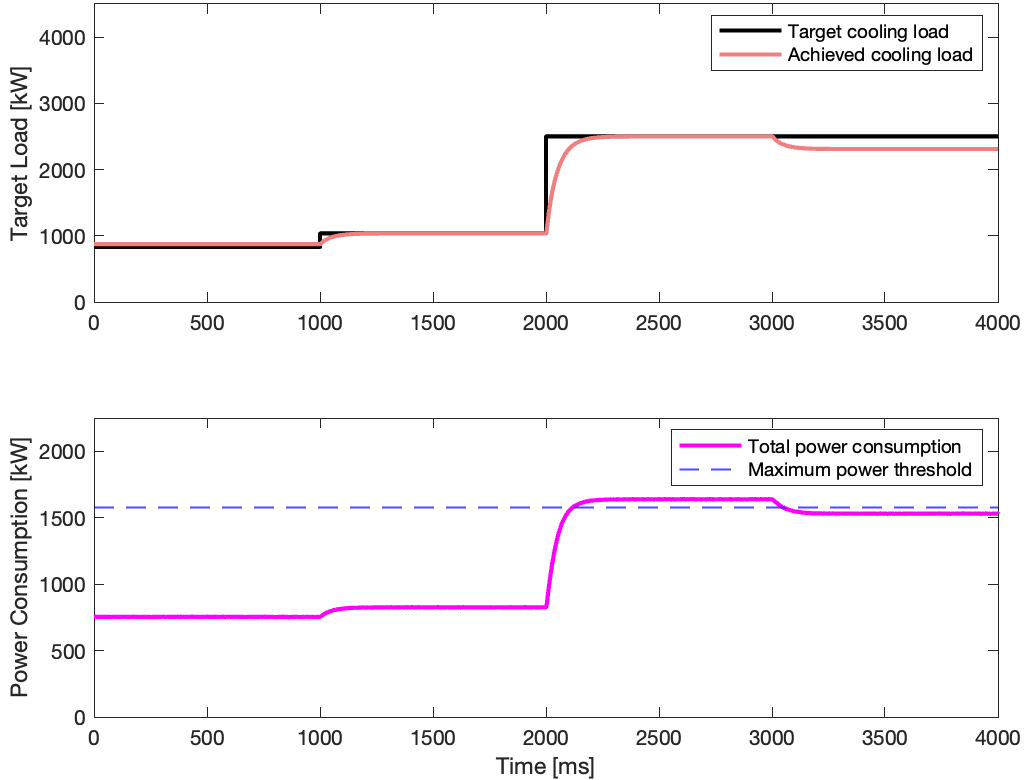}    
    \caption{Safety threshold is exceeded when the system moves far into a previously unexplored space} 
    \label{fig:Figure8}
    \end{center}
\end{figure}

Fig. \ref{Figure4} shows that at the first solution instance of the algorithm, not all true power consumption curves (in green) for the compressors - particularly for unexplored high load regions - lie inside the confidence bounds constructed by GP regression (blue areas). Since the cooling load scenario in Fig. \ref{Figure5-6} does not move rapidly from a low load to a high load, the algorithm gets the chance to update the GPs and the corresponding bounds before a high load is requested. We now investigate a counterexample where the system receives a high cooling load demand earlier, without exploring intermediate loads. In this case, ARTEO moves the compressors to an operating point past the allowed maximum power consumption to track the requested cooling load based on the satisfaction of the UCB-based constraint but steps back from the target load to satisfy the actual safety limit once the steady-state detection algorithm registers the new data points and updates the GPs and the UCBs. Fig. \ref{fig:Figure8} illustrates this counterexample case. A basic approach to avoid the safety violation in this example can be to constrain the search space to the vicinity of the previously explored operating points using a distance measure, which would act analogous to input rate-of-change constraints in optimal control applications.

\subsection{Hyperparameter selection}
\label{sec:HyperparameterSelection}
The performance of ARTEO is sensitive to two factors: (i) the estimation of GP models and (ii) the value of $z$ used for exploration. Gaussian process regression is a technique driven by its own hyperparameters such as kernel selection and noise levels. In this case study - based on hyperparameter optimization - the squared-exponential kernel is chosen as the kernel function for all compressors.
\begin{figure}[b]
    \begin{center}
    \includegraphics[scale=0.4]{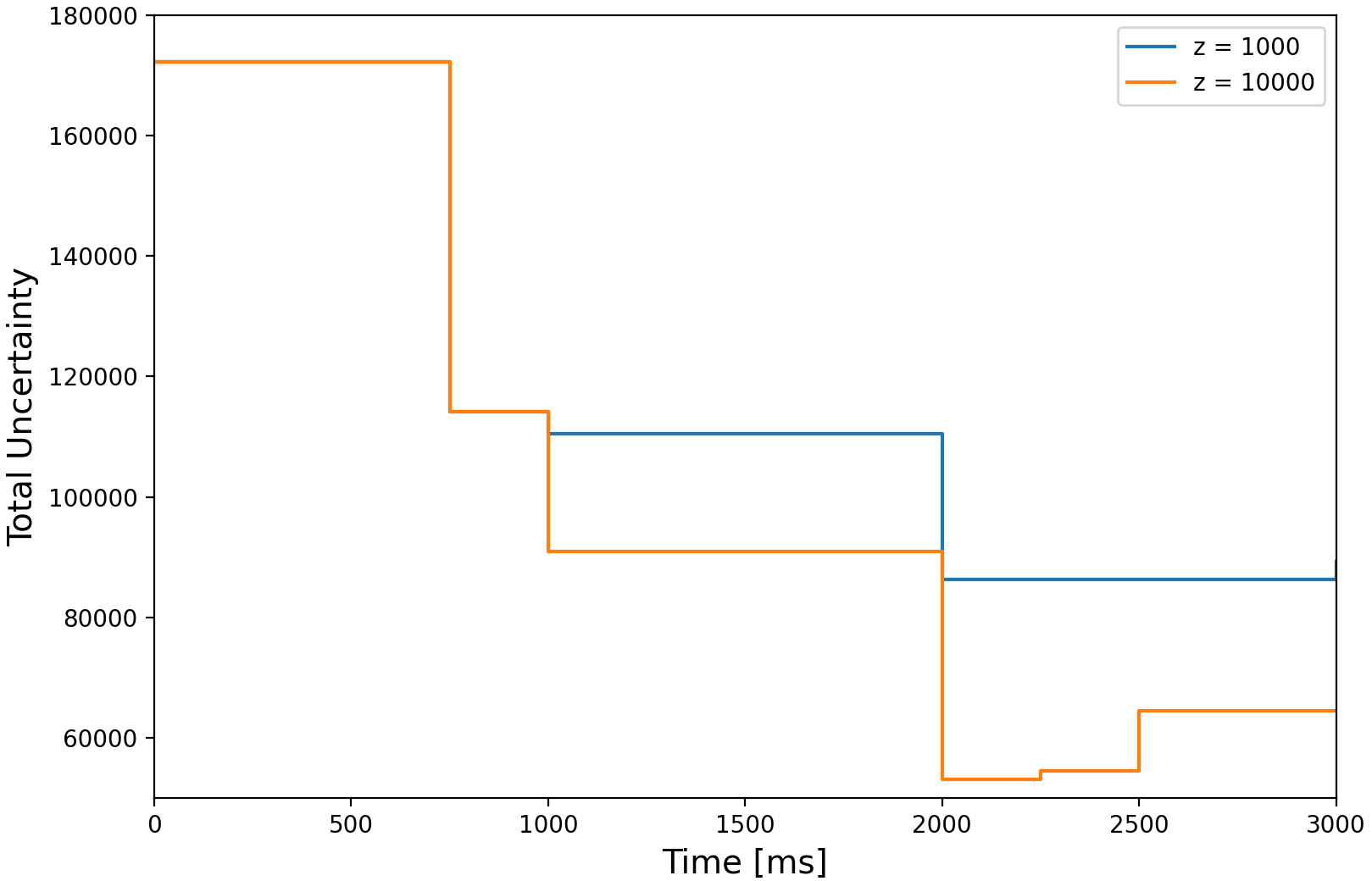}    
    \caption{Total uncertainty decrease in the environment by varying $z$ values} 
    \label{fig:Figure10}
    \end{center}
\end{figure}
We can visualize the total uncertainty in an environment by the standard deviation of GP models summed over a discretization of the operating envelope. In a static environment, we expect the total uncertainty to decrease over time due to exploration. This decrease is led by the size of $z$, with greater z corresponding to a quicker decrease. Fig. \ref{fig:Figure10} shows the trend in the decrease of the total uncertainty in a case study for two different values of $z$. Here, the total cooling demand is kept at a fixed value and individual compressor loads are moved by the ARTEO algorithm to decrease the uncertainty while maintaining the safety threshold. By adjusting the value of $z$, it is possible to shift the emphasis in the utility function to minimisation of the uncertainty, if the term $z$ is greater relative to the power consumption term. Conversely, setting a small value of $z$ enables focusing on the minimisation of power consumption, at the expense of leaving a larger uncertainty.
It is possible to adapt $z$ in a dynamic environment with changing characteristics e.g. when the compressor performances deteriorate over time. This can be achieved by varying the value of $z$ based on the observed prediction errors for the compressor power consumption values. Higher errors will indicate increasing plant-model mismatch and will increase z for more exploration. Forgetting older data points will also increase uncertainty and will drive exploration without changing the value of z.

\begin{figure}[t]
    \begin{center}
    \includegraphics[width=8.4cm]{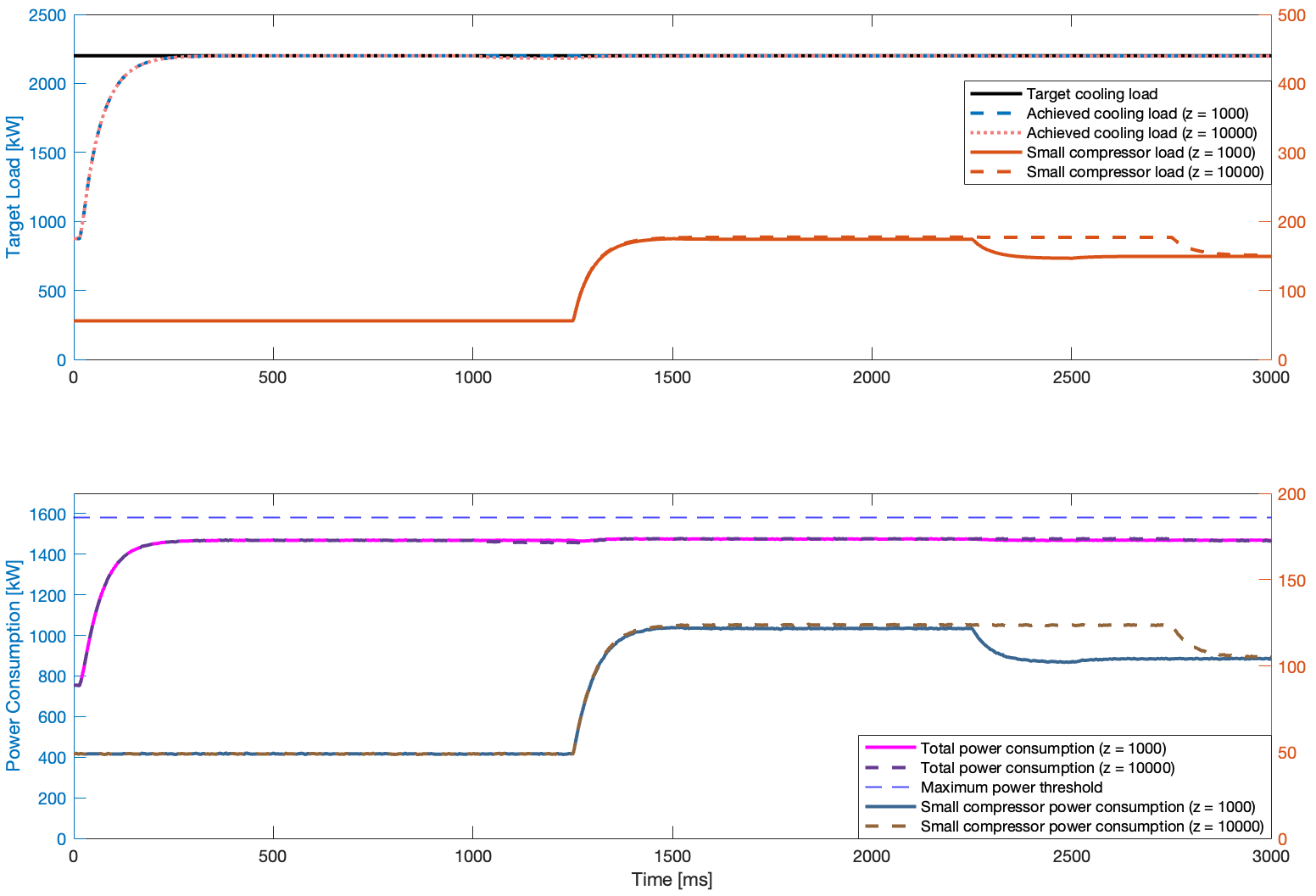}    
    \caption{Comparison of power consumption and target satisfaction for different $z$ values} 
    \label{fig:Figure11}
    \end{center}
\end{figure}

\subsection{Comparison with exact knowledge}

Lastly, we compare ARTEO with an optimization, which has access to exact knowledge of the environment. For this comparison the cooling scenario in Fig. \ref{Figure3} is used again and the optimization problem is solved with the same polynomial equations of the power consumption characteristics used in the simulations. The comparison of the results for power consumption and cooling load demand satisfaction is given in Fig. \ref{Figure5-6}. Here we can see that ARTEO performs close to the optimization with exact knowledge even though the characteristics are initially unknown to ARTEO and are estimated through noisy observations.

\section{Conclusion}
In this work, we introduced and studied an industrial refrigeration case study as an example for the implementation of the ARTEO algorithm for the adaptive and explorative real-time optimization of a safety-critical system under uncertainty. The ARTEO algorithm successfully minimizes the power consumption while tracking the desired cooling load within the operational and safety constraints considered in the main case study. Even though the system starts with very limited information about the process, online learning provides comparable performance to a solution with exact knowledge of the environment.
\bibliography{bibliographyfile.bib}             
                                                   
\end{document}